# THE LAW OF THE SUPREMUM OF A STABLE LÉVY PROCESS WITH NO NEGATIVE JUMPS

By Violetta Bernyk, Robert C. Dalang[1] and Goran Peskir[2]

*Ecole Polytechnique Fédérale, Lausanne, Ecole Polytechnique Fédérale, Lausanne and University of Manchester*

Let $X = (X_t)_{t \geq 0}$ be a stable Lévy process of index $\alpha \in (1,2)$ with no negative jumps and let $S_t = \sup_{0 \leq s \leq t} X_s$ denote its running supremum for $t > 0$. We show that the density function $f_t$ of $S_t$ can be characterized as the unique solution to a weakly singular Volterra integral equation of the first kind or, equivalently, as the unique solution to a first-order Riemann–Liouville fractional differential equation satisfying a boundary condition at zero. This yields an explicit series representation for $f_t$. Recalling the familiar relation between $S_t$ and the first entry time $\tau_x$ of $X$ into $[x, \infty)$, this further translates into an explicit series representation for the density function of $\tau_x$.

**1. Introduction.** In our study [3] of optimal prediction for a stable Lévy process $X = (X_t)_{t \geq 0}$, we encountered the question of computing the distribution function of $S_t = \sup_{0 \leq s \leq t} X_s$ for $t > 0$. In the existing literature, such expressions seem to be available only when $X$ has *no positive jumps* and the purpose of the present paper is to seek similar expressions when $X$ has *no negative jumps*. We note that the latter problem dates back to [5], page 282.

Our main result (Theorem 1) characterizes the density function $f$ of $S_1$ as the unique solution to a weakly singular Volterra integral equation of the first kind or, equivalently, as the unique solution to a first order Riemann–Liouville fractional differential equation satisfying a boundary condition at

Received March 2007; revised October 2007.
[1]Supported by a grant from the Swiss National Foundation for Scientific Research.
[2]Supported by a grant from the British Engineering and Physical Sciences Research Council (EPSRC).

zero. This characterization yields an explicit series representation for $f$ (which, in the case of a Brownian motion, coincides with the well-known expression arising from the reflection principle).

Using the scaling property of $X$, the result extends to $S_t$ for $t \neq 1$. Recalling the familiar relation between $S_t$ and the first entry time $\tau_x$ of $X$ into $[x, \infty)$, this further translates into an explicit series representation for the density function of $\tau_x$ for $x > 0$. Moreover, using the Laplace inversion formula, we derive an integral representation for $f$ (Corollary 2). Finally, we note (Corollary 3) that the proof yields exact constants in the known asymptotic expressions for $f$ at zero and infinity. The knowledge of these constants plays a key role in our treatment of the optimal prediction problem [3].

**2. The result and proof.** 1. Let $X = (X_t)_{t \geq 0}$ be a stable Lévy process of index $\alpha \in (1, 2)$ whose characteristic function is given by

$$(2.1) \quad \mathsf{E}e^{i\lambda X_t} = \exp\left( t \int_0^\infty (e^{i\lambda x} - 1 - i\lambda x) \frac{dx}{\Gamma(-\alpha) x^{1+\alpha}} \right) = e^{t(-i\lambda)^\alpha}$$

for $\lambda \in \mathbb{R}$ and $t \geq 0$. It follows that the Laplace transform of $X$ is given by

$$(2.2) \quad \mathsf{E}e^{-\lambda X_t} = e^{t\lambda^\alpha}$$

for $\lambda \geq 0$ and $t \geq 0$ (the left-hand side being $+\infty$ for $\lambda < 0$). From (2.1) and (2.2), we see that the characteristic exponent of $X$ equals $\Psi(\lambda) = (-i\lambda)^\alpha$, the Laplace exponent of $X$ equals $\psi(\lambda) = \lambda^\alpha$ for $\lambda \geq 0$ and $\varphi(p) := \psi^{-1}(p) = p^{1/\alpha}$ for $p \geq 0$.

2. The following properties of $X$ are readily deduced from (2.1) and (2.2) using standard means (see, e.g., [4] and [13]): the law of $(X_{ct})_{t \geq 0}$ is the same as the law of $(c^{1/\alpha} X_t)_{t \geq 0}$ for each $c > 0$ given and fixed (scaling property); $X$ is a martingale with $\mathsf{E}X_t = 0$ for all $t \geq 0$; $X$ jumps upward (only) and creeps downward [in the sense that $\mathsf{P}(X_{\rho_x} = x) = 1$ for $x < 0$, where $\rho_x = \inf\{t \geq 0 : X_t < x\}$ is the first entry time of $X$ into $(-\infty, x)$]; $X$ has sample paths of unbounded variation; $X$ oscillates from $-\infty$ to $+\infty$ (in the sense that $\liminf_{t \to \infty} X_t = -\infty$ and $\limsup_{t \to \infty} X_t = +\infty$, both a.s.); the starting point $0$ of $X$ is regular [for both $(-\infty, 0)$ and $(0, +\infty)$]. Note that the constant $c = 1/\Gamma(-\alpha)$ in the Lévy measure $\nu(dx) = (c/x^{1+\alpha})\, dx$ of $X$ is chosen/fixed for convenience so that $X$ converges in law to $\sqrt{2}B$ as $\alpha \uparrow 2$, where $B$ is a standard Brownian motion, and all the facts below can be extended to the general constant $c > 0$ depending on $\alpha$ if needed (see Remark 2 below).

3. Let $S_t = \sup_{0 \leq s \leq t} X_s$ denote the running supremum of $X$ for $t \geq 0$ and let $\tau_x = \inf\{t \geq 0 : X_t \geq x\}$ be the first entry time of $X$ into $[x, \infty)$ for $x > 0$. Since $X_s \geq X_{s-}$ for all $s \in [0, t]$ and $X$ is right-continuous, one sees that $\mathsf{P}(S_t \geq x) = \mathsf{P}(\tau_x \leq t)$, so the law of $S_t$ follows from the law of $\tau_x$ (and vice versa). If $X$ is a Lévy process with no positive jumps, then it is known that the two measures

$$\text{(2.3)} \qquad t\mathsf{P}(\tau_x \in dt)\, dx = x\mathsf{P}(X_t \in dx)\, dt$$

coincide on the Borel $\sigma$-algebra of $\mathbb{R}_+ \times \mathbb{R}_+$ (see, e.g., [4], page 190, or [7] and the references therein). This implies that the law of $X_t$ yields the law of $\tau_x$. It follows, in particular, that the known series representations for the density function of $X_t$ (see, e.g., [17], pages 87–89) lead to series representations for the density function of $S_t$. If $X$ has no negative jumps, however, then the identity (2.3) breaks down and no series representation for the density function of $S_t$ seems to be available in the literature. We mention, however, that there is a literature on the distribution of $S_\sigma$ when $\sigma$ is an independent and exponentially distributed random variable, the process $X$ has arbitrary negative jumps, and its positive jumps form a compound Poisson process with the jump-size distribution of the so-called "phase type" (see, e.g., [14] and [2]).

4. Our main result can be stated as follows. Note that $S_t \stackrel{\text{law}}{=} t^{1/\alpha} S_1$ by the scaling property of $X$ so that there is no restriction in assuming that $t = 1$ in the sequel. Recall, also, that $\mathbb{D}^{\alpha-1}$ denotes the *Riemann–Liouville fractional derivative* of order $\alpha - 1$ given by

$$\text{(2.4)} \qquad \mathbb{D}^{\alpha-1} f(x) = \frac{1}{\Gamma(2-\alpha)} \frac{d}{dx} \int_0^x \frac{f(y)}{(x-y)^{\alpha-1}}\, dy$$

for $x > 0$ and any (admissible) function $f : \mathbb{R}_+ \to \mathbb{R}$ (for more details, see, e.g., [16], pages 449–452, and [15], Chapter 2).

THEOREM 1. *Let $X = (X_t)_{t \geq 0}$ be a stable Lévy process of index $\alpha \in (1, 2)$ satisfying (2.1) and (2.2), and let $S_1 = \sup_{0 \leq t \leq 1} X_t$ denote its supremum over the time interval $[0, 1]$. Then the density function $f$ of $S_1$ can be characterized as the unique solution to the weakly singular Volterra integral equation of the first kind*

$$\text{(2.5)} \qquad \int_0^x \left( y + \frac{\alpha}{\Gamma(2-\alpha)} \frac{1}{(x-y)^{\alpha-1}} \right) f(y)\, dy = \frac{\alpha}{\Gamma(1/\alpha)}$$

*or, equivalently, as the unique solution to the fractional differential equation*

$$\text{(2.6)} \qquad xf(x) + \alpha \mathbb{D}^{\alpha-1} f(x) = 0$$

*satisfying the boundary condition*

$$\text{(2.7)} \qquad \lim_{x \downarrow 0} x^{2-\alpha} f(x) = \frac{1}{\Gamma(\alpha-1)\Gamma(1/\alpha)},$$

where $\mathbb{D}^{\alpha-1}$ denotes the Riemann–Liouville fractional derivative given by (2.4) above. This yields the series representation

$$f(x) = \sum_{n=1}^{\infty} \frac{1}{\Gamma(\alpha n - 1)\Gamma(-n + 1 + 1/\alpha)} x^{\alpha n - 2} \tag{2.8}$$

for $x > 0$.

PROOF. To connect the present result with the existing theory, we will begin by recalling a number of known facts about Lévy processes with no positive jumps (for further details, see, e.g., [4], Chapter VII, and [13], Chapter 8).

Let $\widetilde{X} = (\widetilde{X}_t)_{t \geq 0}$ be a Lévy process with no positive jumps starting at zero, let $\widetilde{\Psi}$ denote its characteristic exponent, let $\widetilde{\psi}$ denote the Laplace exponent of $-\widetilde{X}$ and let $\widetilde{\varphi} := \widetilde{\psi}^{(-1)}$ denote the (right) inverse of $\widetilde{\psi}$. Thus, the characteristic function of $\widetilde{X}$ is given by

$$\mathsf{E} e^{i\lambda \widetilde{X}_t} = e^{t\widetilde{\Psi}(\lambda)} \tag{2.9}$$

for $\lambda \in \mathbb{R}$ and the Laplace transform of $-\widetilde{X}$ is given by

$$\mathsf{E} e^{\lambda \widetilde{X}_t} = e^{t\widetilde{\psi}(\lambda)} \tag{2.10}$$

for $\lambda \geq 0$ and $t \geq 0$. Let

$$\widetilde{S}_t = \sup_{0 \leq s \leq t} \widetilde{X}_s \quad \text{and} \quad \widetilde{I}_t = \inf_{0 \leq s \leq t} \widetilde{X}_s \tag{2.11}$$

for $t \geq 0$ and set

$$\widetilde{\tau}_x = \inf\{t \geq 0 : \widetilde{X}_t > x\} \tag{2.12}$$

for $x \geq 0$, on assuming that the infimum is finite a.s.

From the fact that $(e^{\widetilde{\varphi}(p)\widetilde{X}_t - pt})_{t \geq 0}$ is a martingale (and that $\widetilde{X}$ creeps upward), one finds, using the optional sampling theorem, that the Laplace transform of $\widetilde{\tau}_x$ equals

$$\mathsf{E} e^{-p\widetilde{\tau}_x} = e^{-x\widetilde{\varphi}(p)} \tag{2.13}$$

for $p \geq 0$ and $x \geq 0$. Moreover, if $\sigma_p$ is an exponentially distributed random variable with parameter $p > 0$, meaning that $\mathsf{P}(\sigma_p \in dt) = p e^{-pt} dt$ for $t > 0$, which, moreover, is independent of $\widetilde{X}$, then (2.13) implies that

$$\mathsf{P}(\widetilde{S}_{\sigma_p} > x) = \mathsf{P}(\widetilde{\tau}_x \leq \sigma_p) = \mathsf{E} e^{-p\widetilde{\tau}_x} = e^{-x\widetilde{\varphi}(p)} \tag{2.14}$$

for $p > 0$ and $x \geq 0$. This shows that $\widetilde{S}_{\sigma_p}$ is exponentially distributed with parameter $\widetilde{\varphi}(p)$. Hence, one finds that

$$\mathsf{E} e^{\lambda \widetilde{S}_{\sigma_p}} = \frac{\widetilde{\varphi}(p)}{\widetilde{\varphi}(p) - \lambda} \tag{2.15}$$

for $p > 0$ and $\lambda \in \mathbb{C}$ with $\Re(\lambda) < \widetilde{\varphi}(p)$.

Invoking the Wiener–Hopf factorization (see, e.g., [4], page 165, or [13], Theorem 6.16)

$$(2.16) \qquad \mathsf{E}e^{i\lambda \widetilde{X}_{\sigma_p}} = \mathsf{E}e^{i\lambda \widetilde{S}_{\sigma_p}} \mathsf{E}e^{i\lambda \widetilde{I}_{\sigma_p}} = \frac{p}{p - \widetilde{\Psi}(\lambda)},$$

it follows, using (2.15), that

$$(2.17) \qquad \mathsf{E}e^{\lambda \widetilde{I}_{\sigma_p}} = \frac{p(\widetilde{\varphi}(p) - \lambda)}{\widetilde{\varphi}(p)(p - \widetilde{\psi}(\lambda))}$$

for $\lambda \geq 0$ and $p > 0$, on recalling that $\widetilde{\Psi}(-i\lambda) = \widetilde{\psi}(\lambda)$ for $\lambda \geq 0$. The identity (2.17) is well known (see, e.g., [4], page 192, or [13], page 213).

Clearly, $X$ has no negative jumps if and only if $\widetilde{X} := -X$ has no positive jumps, so, by focusing on the left-hand side of (2.17), one finds

$$
\begin{aligned}
\mathsf{E}e^{\lambda \widetilde{I}_{\sigma_p}} &= p \int_0^\infty \mathsf{E}(e^{\lambda \widetilde{I}_t}) e^{-pt}\, dt = p \int_0^\infty \mathsf{E}(e^{-\lambda S_t}) e^{-pt}\, dt \\
(2.18) \qquad &= p \int_0^\infty \left[ 1 - \lambda \int_0^\infty e^{-\lambda x} \mathsf{P}(S_t > x)\, dx \right] e^{-pt}\, dt \\
&= 1 - p\lambda \int_0^\infty e^{-pt}\, dt \int_0^\infty e^{-\lambda x} \mathsf{P}(S_t > x)\, dx
\end{aligned}
$$

for $\lambda \geq 0$ and $p > 0$. Combining (2.17) and (2.18) and noticing/recalling that $\widetilde{\psi}(\lambda) = \psi(\lambda) = \lambda^\alpha$ and $\widetilde{\varphi}(p) = \varphi(p) = p^{1/\alpha}$, one finds that the (joint) time–space Laplace transform of $(t, x) \mapsto \mathsf{P}(S_t > x)$ equals

$$(2.19) \qquad \int_0^\infty e^{-\lambda x}\, dx \int_0^\infty e^{-pt} \mathsf{P}(S_t > x)\, dt = \frac{1}{p - \lambda^\alpha} \left( \frac{1}{p^{1/\alpha}} - \frac{\lambda^{\alpha-1}}{p} \right)$$

for $\lambda > 0$ and $p > 0$.

Note that this formula can also be obtained by taking the Laplace transform with respect to the space variable $x$ on both sides of the expression

$$
\begin{aligned}
\int_0^\infty e^{-pt} \mathsf{P}(S_t > x)\, dt &= \int_0^\infty e^{-pt} \mathsf{P}(\tau_x \leq t)\, dt = \frac{1}{p} \mathsf{E} e^{-p\tau_x} \\
(2.20) \qquad &= \sum_{n=0}^\infty \frac{p^{n-1} x^{\alpha n}}{\Gamma(1 + \alpha n)} - \sum_{n=1}^\infty \frac{p^{n-1-1/\alpha} x^{\alpha n - 1}}{\Gamma(\alpha n)},
\end{aligned}
$$

where the final identity follows from (8.6) in [13], page 214, combined with (ii) and (iii) in [13], page 233. This remark is relevant since the customary approach leading to the closed-form expression (2.20) via the so-called scale function (cf. [13], pages 214–215) corresponds to Laplace inversion (at least formally) with respect to the space parameter. The derivation given below takes a different route by firstly performing Laplace inversion with respect to the time parameter and then dealing with the resulting expression using techniques of linear integral equations (fractional calculus).

After these introductory remarks, we are now ready to move to the first step of the proof, taking (2.19) as the initial point. Below, we will let $\mathbb{L}_p^{-1}$ denote the inverse Laplace transform with respect to the time parameter $p$ and $\mathbb{L}_\lambda^{-1}$ denote the inverse Laplace transform with respect to the space parameter $\lambda$.

1. Considering $p > \lambda^\alpha$ with $\lambda > 0$ fixed, by (3) in [9], page 238, we find

$$(2.21) \qquad \mathbb{L}_p^{-1}\left[\frac{1}{(p-\lambda^\alpha)p^{1/\alpha}}\right](t) = \frac{1}{\Gamma(1/\alpha)}\frac{e^{\lambda^\alpha t}}{\lambda}\gamma(1/\alpha, \lambda^\alpha t)$$

for $t \geq 0$, where $(a,x) \mapsto \gamma(a,x)$ denotes the incomplete gamma function

$$(2.22) \qquad \gamma(a,x) = \int_0^x y^{a-1} e^{-y}\, dy$$

for $a > 0$ and $x \geq 0$. Likewise, by (5) in [9], page 229, we find

$$(2.23) \qquad \mathbb{L}_p^{-1}\left[\frac{\lambda^{\alpha-1}}{(p-\lambda^\alpha)p}\right](t) = \frac{1}{\lambda}(e^{\lambda^\alpha t} - 1)$$

for $t \geq 0$. Combining (2.21) and (2.23), we get

$$(2.24) \qquad \mathbb{L}_p^{-1}\left[\frac{1}{p-\lambda^\alpha}\left(\frac{1}{p^{1/\alpha}} - \frac{\lambda^{\alpha-1}}{p}\right)\right](t) = \frac{1}{\lambda} - \frac{e^{\lambda^\alpha t}}{\lambda}\frac{\Gamma(1/\alpha, \lambda^\alpha t)}{\Gamma(1/\alpha)}$$

for $t \geq 0$, where $(a,x) \mapsto \Gamma(a,x)$ denotes the incomplete gamma function

$$(2.25) \qquad \Gamma(a,x) = \int_x^\infty y^{a-1} e^{-y}\, dy = \Gamma(a) - \gamma(a,x)$$

for $a > 0$ and $x \geq 0$. Since the right-hand side of (2.24) defines a bounded function of $t \geq 0$ and the argument of $\mathbb{L}_p^{-1}$ on the left-hand side is a Laplace transform defined for all $p > 0$ [recall (2.19) above], we see that the identity (2.24) holds globally for $t \geq 0$ and $\lambda > 0$.

2. Note that

$$(2.26) \qquad \begin{aligned}\frac{e^{\lambda^\alpha t}}{\lambda}\frac{\Gamma(1/\alpha, \lambda^\alpha t)}{\Gamma(1/\alpha)} &= \frac{1}{\Gamma(1/\alpha)}\frac{e^{\lambda^\alpha t}}{\lambda}\int_{\lambda^\alpha t}^\infty x^{-1+1/\alpha} e^{-x}\, dx \\ &= \frac{\alpha}{\Gamma(1/\alpha)}\frac{1}{\lambda}e^{t\lambda^\alpha}\int_{t^{1/\alpha}\lambda}^\infty e^{-z^\alpha}\, dz\end{aligned}$$

for $\lambda > 0$ and $t \geq 0$, on substituting $x = z^\alpha$ to obtain the second equality. The final expression in (2.26) reveals a connection with the standard normal distribution corresponding to $\alpha = 2$. Indeed, by the scaling property, it is no restriction to assume that $t = 1$ so that the final expression in (2.26) with $\alpha = 2$ reads

$$\text{(2.27)} \qquad \frac{2}{\sqrt{\pi}} \frac{1}{\lambda} e^{\lambda^2} \int_\lambda^\infty e^{-z^2} \, dz = \frac{e^{\lambda^2}}{\lambda} \text{erfc}(\lambda)$$

for $\lambda > 0$. By (1) in [9], page 265, one knows that

$$\text{(2.28)} \qquad \mathbb{L}_\lambda^{-1}[e^{\lambda^2} \text{erfc}(\lambda)](x) = \frac{1}{\sqrt{\pi}} e^{-x^2/4}$$

and hence it follows that

$$\text{(2.29)} \qquad \mathbb{L}_\lambda^{-1}\left[\frac{e^{\lambda^2}}{\lambda} \text{erfc}(\lambda)\right](x) = \frac{1}{\sqrt{\pi}} \int_0^x e^{-y^2/4} \, dy$$

for $x \geq 0$. The density function $f$ of $S_1$ obtained on the right-hand side of (2.28) and the distribution function $F$ of $S_1$ given on the right-hand side of (2.29) coincide with the expressions obtained from the reflection principle $M_1 := \max_{0 \leq t \leq 1} B_t =^{\text{law}} |B_1|$, which yields $S_1 =^{\text{law}} \sqrt{2} M_1 =^{\text{law}} \sqrt{2} |B_1| =^{\text{law}} |B_2|$, where $B = (B_t)_{t \geq 0}$ is a standard Brownian motion.

3. By the scaling property, it is no restriction to assume that $t = 1$ in the sequel. Let $F$ denote the distribution function of $S_1$ and let $f$ denote the density function of $S_1$. Note that the form of the Laplace transform on the right-hand side of (2.24) being combined with (2.26) implies that the density function exists (see the proof of Corollary 2 below for further detail). Combining (2.19), (2.24) and (2.26), on using that $\mathbb{L}_\lambda^{-1}[\frac{1}{\lambda}\mathbb{L}[f](\lambda)](x) = F(x)$ since $F(x) = \int_0^x f(y) \, dy$ for $x \geq 0$, it follows that

$$\text{(2.30)} \qquad f(x) = \frac{\alpha}{\Gamma(1/\alpha)} \mathbb{L}_\lambda^{-1}\left[e^{\lambda^\alpha} \int_\lambda^\infty e^{-z^\alpha} \, dz\right](x)$$

for $x \geq 0$.

To simplify the notation, consider the equation

$$\text{(2.31)} \qquad g(x) = \mathbb{L}_\lambda^{-1}[G(\lambda)](x)$$

for $x > 0$, where we set

$$\text{(2.32)} \qquad G(\lambda) = e^{\lambda^\alpha} \int_\lambda^\infty e^{-z^\alpha} \, dz$$

for $\lambda > 0$. From (2.32), we see that $G'(\lambda) = \alpha \lambda^{\alpha-1} G(\lambda) - 1$ so that

$$\text{(2.33)} \qquad \frac{G'(\lambda)}{\lambda} - \frac{\alpha}{\lambda^{2-\alpha}} G(\lambda) + \frac{1}{\lambda} = 0$$

for $\lambda > 0$. Since $\mathbb{L}_\lambda^{-1}[G'(\lambda)](x) = -xg(x)$, it follows that $\mathbb{L}_\lambda^{-1}[G'(\lambda)/\lambda](x) = -\int_0^x yg(y)dy$ for $x > 0$. Moreover, using (1) in [9], page 137, we see that $\mathbb{L}_\lambda^{-1}[1/\lambda^{2-\alpha}](x) = 1/(\Gamma(2-\alpha)x^{\alpha-1})$ so that $\mathbb{L}_\lambda^{-1}[G(\lambda)/\lambda^{2-\alpha}](x) = (1/\Gamma(2-$ $\alpha$)) $\int_0^x (g(y)/(x-y)^{\alpha-1})\,dy$ for $x > 0$. Finally, we have $\mathbb{L}_\lambda^{-1}[1/\lambda](x) = 1$ for $x > 0$. Hence, taking $\mathbb{L}_\lambda^{-1}$ in (2.33), we find that

$$-\int_0^x y g(y)\,dy - \frac{\alpha}{\Gamma(2-\alpha)} \int_0^x \frac{g(y)}{(x-y)^{\alpha-1}}\,dy + 1 = 0 \tag{2.34}$$

for $x > 0$. Noting that $g(x) = (\Gamma(1/\alpha)/\alpha)f(x)$ for $x > 0$, we see that (2.34) reads

$$\int_0^x y f(y)\,dy + \frac{\alpha}{\Gamma(2-\alpha)} \int_0^x \frac{f(y)}{(x-y)^{\alpha-1}}\,dy = \frac{\alpha}{\Gamma(1/\alpha)} \tag{2.35}$$

for $x > 0$ and this is exactly equation (2.5).

4. We will seek a solution to (2.35) of the form

$$f(x) = \sum_{n=0}^{\infty} a_n x^{\beta n + \gamma}, \tag{2.36}$$

where $\beta$ and $\gamma$ are constants to be determined. First, note that

$$\int_0^x y f(y)\,dy = \sum_{n=0}^{\infty} a_n \int_0^x y^{\beta n + \gamma + 1}\,dy = \sum_{n=0}^{\infty} \frac{a_n}{\beta n + \gamma + 2} x^{\beta n + \gamma + 2} \tag{2.37}$$

for $x > 0$. Second, by (3.191) in [10], page 333, and (6.2.2) in [1], page 258, we have

$$\int_0^x y^{\mu-1}(x-y)^{\nu-1}\,dy = x^{\mu+\nu-1} B(\mu, \nu) = x^{\mu+\nu-1} \frac{\Gamma(\mu)\Gamma(\nu)}{\Gamma(\mu+\nu)} \tag{2.38}$$

for $\mu > 0$, $\nu > 0$ and $x > 0$. It follows that

$$\begin{aligned}
\frac{\alpha}{\Gamma(2-\alpha)} \int_0^x \frac{f(y)}{(x-y)^{\alpha-1}}\,dy &= \frac{\alpha}{\Gamma(2-\alpha)} \sum_{n=0}^{\infty} a_n \int_0^x \frac{y^{\beta n + \gamma}}{(x-y)^{\alpha-1}}\,dy \\
&= \sum_{n=0}^{\infty} a_n \frac{\alpha \Gamma(\beta n + \gamma + 1)}{\Gamma(\beta n + \gamma - \alpha + 3)} x^{\beta n + \gamma - \alpha + 2}
\end{aligned} \tag{2.39}$$

for $x > 0$. Combining (2.35), (2.37) and (2.39), we find that $\beta = \alpha$ and $\gamma = \alpha - 2$. Inserting (2.37) and (2.39) into (2.35) with these $\beta$ and $\gamma$, we get

$$\sum_{n=0}^{\infty} (a_n A_n + a_{n+1} B_{n+1}) x^{\alpha(n+1)} + a_0 B_0 = \frac{\alpha}{\Gamma(1/\alpha)}, \tag{2.40}$$

where the constants $A_n$ and $B_n$ are defined by

$$A_n = \frac{1}{\alpha(n+1)} \quad \text{and} \quad B_n = \alpha \frac{\Gamma(\alpha(n+1)-1)}{\Gamma(\alpha n + 1)} \tag{2.41}$$

for $n \geq 0$. From (2.40) and (2.41), we find, by induction, that

$$\text{(2.42)} \qquad a_n = (-1)^n \frac{A_{n-1} A_{n-2} \cdots A_1 A_0}{B_n B_{n-1} \cdots B_2 B_1} a_0$$

for $n \geq 1$, where $a_0 = 1/(\Gamma(1/\alpha)\Gamma(\alpha - 1))$. Inserting (2.42) into (2.36) with $\beta = \alpha$ and $\gamma = \alpha - 2$, and making use of (2.41), we obtain the series representation

$$\text{(2.43)} \qquad \begin{aligned} f(x) = & \frac{1}{\Gamma(1/\alpha)} \\ & \times \sum_{n=0}^{\infty} \frac{(-1)^n}{\alpha^{2n} n!} \frac{\Gamma(n\alpha + 1)\Gamma((n-1)\alpha + 1) \cdots \Gamma(\alpha + 1)\Gamma(1)}{\Gamma((n+1)\alpha - 1)\Gamma(n\alpha - 1) \cdots \Gamma(2\alpha - 1)\Gamma(\alpha - 1)} \\ & \times x^{\alpha(n+1) - 2} \end{aligned}$$

for $x > 0$. Using Stirling's formula $\Gamma(ax + b) \sim \sqrt{2\pi} e^{-ax}(ax)^{ax+b-1/2}$ as $x \to \infty$, where $a > 0$ and $b \in \mathbb{R}$ (cf. (6.1.39) in [1], page 257), it is readily verified that $|a_{n+1}/a_n| = O(n^{1-\alpha})$ as $n \to \infty$, whence the ratio test implies that the series in (2.43) converges absolutely for every $x > 0$ and that $f$ defined by (2.43) is a continuous function on $(0, \infty)$. Note, also, that only the leading term $(1/(\Gamma(\alpha - 1)\Gamma(1/\alpha))x^{\alpha-2}$ of the series is singular at zero, so we can integrate in (2.43) term by term over any finite interval in $[0, \infty)$. Finally, by induction over $n \geq 0$, using the fact that $\Gamma(x + 1) = x\Gamma(x)$ for $x \in \mathbb{R} \setminus \{0, -1, -2, \ldots\}$, it is easily verified that the series representation (2.43) can be simplified to the form given in (2.8) above.

5. We now show that $f$ from (2.8) is a unique solution to the integral equation (2.5). For this, let us first note that since $f$ satisfies (2.43) and hence solves (2.35), it follows that $g = (\Gamma(1/\alpha)/\alpha)f$ solves (2.34). Assuming that $g$ has a Laplace transform and taking the Laplace transform $\mathbb{L}$ on both sides of (2.34), we see that $G = \mathbb{L}[g]$ solves (2.33). The general solution to (2.33) is given by $G(\lambda) = ce^{\lambda^\alpha} + e^{\lambda^\alpha} \int_\lambda^\infty e^{-z^\alpha} dz$ for $\lambda > 0$, where $c$ is a constant. In order to compute the Laplace transform of $f$ defined in (2.43), we could attempt to interchange $\mathbb{L}$ and the sum and use the fact that $\mathbb{L}[x^\rho](\lambda) = \Gamma(\rho + 1)/\lambda^{\rho+1}$ for $\rho > -1$ and $\lambda > 0$. Using the ratio test, however, it is possible to verify that the resulting series diverges and therefore is not equal to $\mathbb{L}[f]$. We note, however, that if we could show that $\int_0^\infty e^{-\lambda x} f(x) dx \to 0$ as $\lambda \to \infty$, then we would have $c = 0$ and (2.30) would imply that $f$ from (2.43) is indeed the density function of $S_1$, as claimed.

Given this difficulty, we shall take a different tack and establish uniqueness of the solution to (2.35) in the class of functions that are locally integrable on $[0, \infty)$ and bounded on compact subsets of $(0, \infty)$ [these conditions are natural requirements so that the left-hand side of (2.35) makes sense]. Multiplying both sides of (2.35) by $(z - x)^{\alpha-2}$ and integrating the resulting identity with respect to $x$ from $0$ to $z$, we can use Fubini's theorem and (2.38) to obtain

$$
\begin{aligned}
&\frac{1}{\alpha-1}\int_0^z y(z-y)^{\alpha-1}f(y)\,dy + \alpha\Gamma(\alpha-1)\int_0^z f(y)\,dy \\
&= \frac{\alpha}{(\alpha-1)\Gamma(1/\alpha)}z^{\alpha-1}
\end{aligned}
\tag{2.44}
$$

for $z > 0$. Note that the interchange of the order of integration above is justified whenever $f$ is locally integrable on $[0,\infty)$ and bounded on compact subsets of $(0,\infty)$. Differentiating this identity with respect to $z$ and substituting $x$ for $z$, we get

$$
\frac{1}{\alpha\Gamma(\alpha-1)}\int_0^x \frac{y}{(x-y)^{2-\alpha}}f(y)\,dy + f(x) = \frac{1}{\Gamma(1/\alpha)\Gamma(\alpha-1)}\frac{1}{x^{2-\alpha}}
\tag{2.45}
$$

for $x > 0$. This is a weakly singular Volterra integral equation of the second kind. Previous considerations show that both the function from (2.43) and the density function from (2.30) solve the equation (2.45). We note in passing that when $f$ is the density function, we see from (2.45) that $f(x) \leq [1/(\Gamma(\alpha-1)\Gamma(1/\alpha))](1/x^{2-\alpha})$ for all $x > 0$ so that $f$ is bounded on compact subsets of $(0,\infty)$.

Denote by $\phi$ the difference between the two solutions to (2.45). Then

$$
a\int_0^x \frac{y}{(x-y)^{2-\alpha}}\phi(y)\,dy + \phi(x) = 0
\tag{2.46}
$$

for $x > 0$, where we set $a = 1/(\alpha\Gamma(\alpha-1))$. It follows from [12], Theorem 7, page 35, that $\phi = 0$ if $\phi$ is locally square-integrable, but since the latter could not be the case (around zero), we give a direct proof of the former fact. For this, fix $x_1 > 0$ arbitrarily large and set $\xi = |\phi|$. Letting

$$
\mathrm{T}\xi(x) = \int_0^x \frac{\xi(y)}{(x-y)^{2-\alpha}}\,dy,
\tag{2.47}
$$

we find, by induction using (2.46), that

$$
\xi(x) \leq b^n \mathrm{T}^n \xi(x)
\tag{2.48}
$$

for $x \in (0, x_1]$ and $n \geq 1$, where $b = ax_1$. An iterative calculation using Fubini's theorem and (2.38) shows that

$$
\mathrm{T}^n \xi(x) \leq c_n \int_0^x \frac{\xi(y)}{(x-y)^{1-n(\alpha-1)}}\,dy
\tag{2.49}
$$

for $x \in (0, x_1]$ and $n \geq 1$ with some constant $c_n > 0$. Choosing $n \geq 1$ large enough so that $1 - n(\alpha-1) < 0$, combining (2.48) with (2.49) and applying a simple iteration procedure to the resulting inequality, we find that

$$
\xi(x) \leq \frac{c^m x^{m-1}}{(m-1)!}\int_0^x \xi(y)\,dy
\tag{2.50}
$$

for $x \in (0, x_1]$ and $m \geq 1$ with some constant $c > 0$. Since the right-hand side converges to zero as $m \to \infty$, it follows that $\xi(x) = 0$ for $x \in (0, x_1]$ and thus $\phi(x) = 0$ for all $x > 0$. This shows that the two solutions to (2.45) coincide on $(0, \infty)$. Hence, we can conclude that $f$ from (2.8) is a unique solution to (2.5) in the class of functions which are locally integrable on $[0, \infty)$ and bounded on compact subsets of $(0, \infty)$.

6. Note that the fractional differential equation (2.6) follows from (2.5) by differentiation so that (2.8) defines its solution satisfying the boundary condition (2.7). Now, suppose that $f$ solves (2.6) and satisfies (2.7). Then (2.5) follows from (2.6) by integration [on using (2.38) with $\mu = \alpha - 1$], so $f$ solves (2.35). Then proceeding as above, we find that $f$ must be equal to the density function, as long as $f$ is locally integrable on $[0, \infty)$ and bounded on compact subsets of $(0, \infty)$. This establishes the existence and uniqueness claim about (2.6) and (2.7) in the latter class of functions. The proof of the theorem is complete. □

REMARK 1. The integral equation (2.5) is closely related to the (generalized) *Abel equation of the first kind*

$$(2.51) \qquad \int_0^x \left( a + \frac{1}{(x-y)^\beta} \right) f(y) \, dy = R(x) \qquad (0 < \beta < 1),$$

which admits a closed-form solution expressed in terms of the Riemann–Liouville fractional derivative of $R$ (of order $1 - \beta$). For more details, see [16] and the references therein. Note that the integral equation (2.5) is of the form

$$(2.52) \qquad \int_0^x \left( ay + \frac{1}{(x-y)^\beta} \right) f(y) \, dy = R(x) \qquad (0 < \beta < 1),$$

which may be viewed as being of the *first order* if the Abel equation (2.51) is viewed as being of the *zeroth order*. Note also that the equation (2.45) is the "second kind" analog of the equation (2.5).

REMARK 2. The results of Theorem 1 extend to the case when the Lévy measure equals

$$(2.53) \qquad \nu(dx) = \frac{c}{x^{1+\alpha}} \, dx,$$

where $c > 0$ is a general constant. This can be derived using the scaling property of $X$. Letting, in this case, $f_t$ denote the density function of $S_t = \sup_{0 \leq s \leq t} X_s$, we note for future reference that (2.8) extends as follows:

$$(2.54) \quad f_t(x) = \sum_{n=1}^\infty \frac{1}{(c\Gamma(-\alpha)t)^{n-1/\alpha} \Gamma(\alpha n - 1) \Gamma(-n + 1 + 1/\alpha)} x^{\alpha n - 2}$$

for $x > 0$ and $t > 0$. Similarly, from (2.30), it is readily verified that

$$\mathsf{E}e^{-\lambda S_t} = \int_0^\infty e^{-\lambda x} f_t(x)\,dx = \frac{\alpha}{\Gamma(1/\alpha)} e^{\kappa t \lambda^\alpha} \int_{(\kappa t)^{1/\alpha}\lambda}^\infty e^{-z^\alpha}\,dz \tag{2.55}$$

for $\lambda > 0$ and $t > 0$, where we set $\kappa = c\Gamma(-\alpha)$. Note, in particular, that (2.55) yields $\mathsf{E}S_t = (\alpha/\Gamma(1/\alpha))(\kappa t)^{1/\alpha}$ for $t > 0$.

7. Further to the series representation given in (2.8) above, the next corollary presents an integral representation for the density function $f$ of $S_1$. Since this representation extends to $t \neq 1$ and $c \neq 1/\Gamma(-\alpha)$ by the scaling property of $X$, we will only focus on the case when $t = 1$ and $c = 1/\Gamma(-\alpha)$ in (2.53). We refer to [11], Theorem 1, page 422, and [6], Theorem 3, page 74, for more general integral representations in this context (with no obvious connection to the one given below).

COROLLARY 2. *Let $X = (X_t)_{t \geq 0}$ be a stable Lévy process of index $\alpha \in (1,2)$ satisfying (2.1) and (2.2), and let $S_1 = \sup_{0 \leq t \leq 1} X_t$ denote its supremum over the time interval $[0,1]$. Then the density function $f$ of $S_1$ is given by*

$$\begin{aligned} f(x) = \frac{1}{\pi} \int_0^\infty \bigg[& e^{t^\alpha \cos(\alpha\pi/2)} \cos(t^\alpha \sin(\alpha\pi/2) + tx) \\ &+ \frac{1}{\Gamma(1/\alpha)} \int_0^{t^\alpha} \frac{e^{y \cos(\alpha\pi/2)}}{(t^\alpha - y)^{1-1/\alpha}} \\ &\qquad \times \sin(y \sin(\alpha\pi/2) + tx)\,dy \bigg] dt \end{aligned} \tag{2.56}$$

*for $x > 0$.*

PROOF. Setting in the first equality and noting in the second equality that

$$H(\lambda) = \int_\lambda^\infty e^{-y^\alpha}\,dy = \frac{\Gamma(1/\alpha)}{\alpha} - \int_0^\lambda e^{-y^\alpha}\,dy \tag{2.57}$$

for $\lambda > 0$, it follows from (2.30) that

$$f(x) = \frac{\alpha}{\Gamma(1/\alpha)} \mathbb{L}_\lambda^{-1}[e^{\lambda^\alpha} H(\lambda)](x) \tag{2.58}$$

for $x > 0$. From the second equality in (2.57), one sees that $H$ can be analytically continued to the entire complex plane. The same fact is therefore true for $\lambda \mapsto e^{\lambda^\alpha} H(\lambda)$ so that the Laplace inversion formula is applicable in

(2.58) yielding

$$
\begin{aligned}
(2.59) \quad f(x) &= \frac{\alpha}{\Gamma(1/\alpha)} \frac{1}{2\pi i} \int_{-i\infty}^{+i\infty} e^{z^\alpha} \left( \frac{\Gamma(1/\alpha)}{\alpha} - \int_0^z e^{-y^\alpha} \, dy \right) e^{xz} \, dz \\
&= \frac{\alpha}{\Gamma(1/\alpha)} \frac{1}{2\pi} \int_{-\infty}^{+\infty} e^{(it)^\alpha} \left( \frac{\Gamma(1/\alpha)}{\alpha} - it \int_0^1 e^{-(ity)^\alpha} \, dy \right) e^{itx} \, dt
\end{aligned}
$$

for $x > 0$. In the case $t > 0$, we have $\exp((it)^\alpha) = \exp(t^\alpha(\cos(\alpha\pi/2) + i\sin(\alpha\pi/2)))$ and $\exp(-(ity)^\alpha) = \exp(-(ty)^\alpha(\cos(\alpha\pi/2) + i\sin(\alpha\pi/2)))$. In the case $t < 0$, we have $\exp((it)^\alpha) = \exp((-t)^\alpha(\cos(\alpha\pi/2) - i\sin(\alpha\pi/2)))$ and $\exp(-(ity)^\alpha) = \exp(-(-ty)^\alpha(\cos(\alpha\pi/2) - i\sin(\alpha\pi/2)))$. Inserting these expressions into (2.59), one can verify that the integral from $-\infty$ to $+\infty$ equals twice the integral from $0$ to $+\infty$, which, in turn, can be reduced to the form given in (2.56) above. As this verification is somewhat lengthy, but still straightforward, further details will be omitted. This completes the proof. □

8. The next corollary describes asymptotic behavior of the law of $S_1$ at zero and infinity. Recall that $f(x) \sim g(x)$ as $x \to x_0$ means that $\lim_{x \to x_0} f(x)/g(x) = 1$ for $x_0 \in [-\infty, +\infty]$.

COROLLARY 3. *Let $X = (X_t)_{t \geq 0}$ be a stable Lévy process of index $\alpha \in (1, 2)$ whose characteristic function is given by*

$$
(2.60) \quad \mathsf{E} e^{i\lambda X_t} = \exp\left( t \int_0^\infty (e^{i\lambda x} - 1 - i\lambda x) \frac{c \, dx}{x^{1+\alpha}} \right)
$$

*for $\lambda \in \mathbb{R}$ and $t \geq 0$, where $c > 0$ is a given and fixed constant. Let $S_1 = \sup_{0 \leq t \leq 1} X_t$ and let $f$ denote the density function of $S_1$. Then*

$$
(2.61) \quad f(x) \sim \frac{1}{(c\Gamma(-\alpha))^{1-1/\alpha} \Gamma(\alpha - 1) \Gamma(1/\alpha)} x^{\alpha-2} \quad \text{as } x \downarrow 0,
$$

$$
(2.62) \quad f(x) \sim c x^{-\alpha-1} \quad \text{as } x \uparrow \infty.
$$

PROOF. The relation (2.61) follows directly from the explicit series representation (2.54). The relation (2.62) can be derived from the integral representation (2.56), as shown in [8]. □

## REFERENCES

[1] ABRAMOWITZ, M. and STEGUN, I. A. (1992). *Handbook of Mathematical Functions.* Dover, New York. MR1225604
[2] ASMUSSEN, S., AVRAM, F. and PISTORIUS, M. R. (2004). Russian and American put options under exponential phase-type Lévy models. *Stochastic Process. Appl.* **109** 79–111. MR2024845

[3] BERNYK, V., DALANG, R. C. and PESKIR, G. (2007). Predicting the ultimate supremum of a stable Lévy process. Research Report No. 28, Probab. Statist. Group Manchester.
[4] BERTOIN, J. (1996). *Lévy Processes.* Cambridge Univ. Press. [MR1406564](MR1406564)
[5] BINGHAM, N. H. (1973). Maxima of sums of random variables and suprema of stable processes. *Z. Wahrsch. Verw. Gebiete* **26** 273–296. [MR0415780](MR0415780)
[6] BOROVKOV, A. A. (1976). *Stochastic Processes in Queueing Theory.* Springer, New York. [MR0391297](MR0391297)
[7] BOROVKOV, K. and BURQ, Z. (2001). Kendall's identity for the first crossing time revisited. *Electron. Comm. Probab.* **6** 91–94. [MR1871697](MR1871697)
[8] DONEY, R. A. (2008). A note on the supremum of a stable process. *Stochastics* **80** 151–155. [MR2402160](MR2402160)
[9] ERDÉLYI, A. (1954). *Tables of Integral Transforms* **1**. McGraw-Hill, New York. [MR0061695](MR0061695)
[10] GRADSHTEYN, I. S. and RYZHIK, I. M. (1994). *Table of Integrals, Series, and Products.* Academic Press, New York. [MR1243179](MR1243179)
[11] HEYDE, C. C. (1969). On the maximum of sums of random variables and the supremum functional for stable processes. *J. Appl. Probab.* **6** 419–429. [MR0251766](MR0251766)
[12] HOCHSTADT, H. (1973). *Integral Equations.* Wiley, New York. [MR0390680](MR0390680)
[13] KYPRIANOU, A. E. (2006). *Introductory Lectures on Fluctuations of Lévy Processes with Applications.* Springer, Berlin. [MR2250061](MR2250061)
[14] MORDECKI, E. (2002). The distribution of the maximum of the Lévy process with positive jumps of phase-type. *Theory Stoch. Process.* **8** 309–316. [MR2027403](MR2027403)
[15] PODLUBNY, I. (1999). *Fractional Differential Equations.* Academic Press, New York. [MR1658022](MR1658022)
[16] POLYANIN, A. D. and MANZHIROV, A. V. (1998). *Handbook of Integral Equations.* Chapman and Hall, Boca Raton, FL. [MR1790925](MR1790925)
[17] SATO, K. (1999). *Lévy Processes and Infinitely Divisible Distributions.* Cambridge Univ. Press. [MR1739520](MR1739520)

V. BERNYK
STATISTICAL LABORATORY
UNIVERSITY OF CAMBRIDGE
WILBERFORCE ROAD
CAMBRIDGE CB3 0WB
UNITED KINGDOM
E-MAIL: V.Bernyk@statslab.cam.ac.uk

R. C. DALANG
INSTITUT DE MATHÉMATIQUES
ECOLE POLYTECHNIQUE FÉDÉRALE
STATION 8
1015 LAUSANNE
SWITZERLAND
E-MAIL: robert.dalang@epfl.ch

G. PESKIR
SCHOOL OF MATHEMATICS
THE UNIVERSITY OF MANCHESTER
OXFORD ROAD
MANCHESTER M13 9PL
UNITED KINGDOM
E-MAIL: goran@maths.man.ac.uk